\documentclass{article}

\usepackage{graphicx}
\usepackage{ulem}
\usepackage{amsmath, amsthm,pgf,tikz}
\usetikzlibrary{arrows}
\usepackage{amsfonts}
\usepackage{amssymb}
\usepackage{url}
\usepackage{hyperref}  %This only works for pdflatex, not latex
\usepackage{verbatim}
\usepackage{color}
\definecolor{red}{rgb}{1,0,0}

\definecolor{blu}{rgb}{0,0,1}

\usepackage[mathlines]{lineno}

\usepackage{pgf,tikz}
\usepackage{mathrsfs}
\usetikzlibrary{arrows}
\definecolor{qqqqff}{rgb}{0.,0.,1.}

\def\noi{\noindent}

\newtheorem{thm}{Theorem}
\newtheorem{cor}[thm]{Corollary}

\theoremstyle{definition}
\newtheorem{rem}[thm]{Remark}

\theoremstyle{definition}

\theoremstyle{definition}

\theoremstyle{definition}

\newcommand{\bit}{\begin{itemize}}
\newcommand{\eit}{\end{itemize}}
\newcommand{\ben}{\begin{enumerate}}
\newcommand{\een}{\end{enumerate}}
\newcommand{\beq}{\begin{equation}}
\newcommand{\eeq}{\end{equation}}
\newcommand{\bea}{\begin{eqnarray*}}
\newcommand{\eea}{\end{eqnarray*}}
\newcommand{\bpf}{\begin{proof}}
\newcommand{\epf}{\end{proof}}

\date{July 9, 2018 }

\begin{document}

\title{Chorded pancyclicity in $k$-partite graphs. }

\author{ Daniela Ferrero\thanks{Department of Mathematics, Texas State University, San Marcos, TX, USA (dferrero@txstate.edu)},
Linda Lesniak\thanks{Department of Mathematics, Western Michigan University, Kalamazoo, MI, USA (linda.lesniak@wmich.edu)}}

\maketitle

\begin{abstract}
 We prove that for any integers $p\geq k\geq 3$ and any $k$-tuple of positive integers $(n_1,\ldots ,n_k)$ such that $p=\sum _{i=1}^k{n_i}$ and $n_1\geq n_2\geq \ldots \geq n_k$, the condition $n_1\leq {p\over 2}$ is necessary and sufficient for every subgraph of the complete $k$-partite graph $K(n_1,\ldots ,n_k)$ with at least \[{{4 -2p+2n_1+\sum _{i=1}^{k}  n_i(p-n_i)}\over 2}\] edges to be chorded pancyclic. Removing all but one edge incident with any vertex of minimum degree in $K(n_1,\ldots ,n_k)$ shows that this result is best possible. Our result implies that for any integers, $k\geq 3$ and $n\geq 1$, a balanced $k$-partite graph of order $kn$ with has at least ${{(k^2-k)n^2-2n(k-1)+4}\over 2}$ edges is chorded pancyclic. In the case $k=3$, this result strengthens a previous one by Adamus, who in 2009 showed that a balanced tripartite graph of order $3n$, $n \geq 2$, with at least $3n^2 - 2n + 2$ edges is pancyclic.
\end{abstract}

\noi {\bf Keywords:} hamiltonicity;  pancyclicity; bipancyclicity; chorded pancycliclity; bipartite graphs; $k$-partite graphs.  

\smallskip

\noi{\bf AMS subject classification} 05C45 

\color{black}
\section{Background}

A graph $G$ is {\it hamiltonian}  if it has a spanning cycle. One of the earliest sufficient conditions for a graph to be hamiltonian is one due to Ore \cite{O60}.

\begin{thm} [Ore 1960] \label{Ore}
	Let $G$ be a graph of order $p\geq 3$. If for every pair $u$ and $v$ of nonadjacent vertices, $d(u)+d(v)\geq p$, then $G$ is hamiltonian.
\end{thm}

Two immediate corollaries of Ore's theorem are a  minimum degree condition  due to Dirac \cite{D52} and a simple edge condition.

\begin{cor} [Dirac 1952]\label{Dirac}
	Let $G$ be a graph of order $p\geq 3$. If $d(v) \geq p/2$ for every vertex $v$ of $G$, then $G$ is hamiltonian.
	
\end{cor} 

\begin{cor} \label{edgecondition1}
	Let $G$ be a graph of order $p \geq 3$.  If $G$ has at least  $\dfrac{p^2-3p+6}{2}$  edges, then $G$ is hamiltonian.
\end{cor}

In \cite{B71} Bondy introduced the notion of pancylicity in graphs.  A graph $G$ of order $p \geq 3$ is {\it pancyclic} if it not only has a spanning cycle as do hamiltonian graphs, but also a cycle of order $t$ for every $3 \leq t \leq p$.  Thus every pancyclic graph is hamiltonian but not necessarily the converse.  However, Bondy showed that our three sufficient conditions for hamiltonicity were ``almost" sufficient for pancyclicity.  

\begin{thm} [Bondy 1971]
	Let $G$ be a graph of order $p\geq 3$. If for every pair $u$ and $v$ of nonadjacent vertices, $d(u)+d(v)\geq p$, then either $G$ is pancyclic or $p$ is even and  $G$ is the complete bipartite graph $K_{p/2,p/2}$.	
\end{thm}

\begin{cor}
		Let $G$ be a graph of order $p\geq 3$. If $d(v) \geq p/2$ for every vertex $v$ of $G$, then either $G$ is pancyclic or $p$ is even and  $G$ is the complete bipartite graph $K_{{p/2},{p/2}}$.	
\end{cor}

\begin{cor}\label{generalpancyclic}
		Let $G$ be a graph of order $p \geq 3$.  If $G$ has at least  $\dfrac{p^2-3p+6}{2}$  edges, then $G$ is pancyclic.
\end{cor}
	
The following result of Bondy \cite{B71} gives a sufficient condition for a hamiltonian graph to be pancyclic that we will refer to later in this paper.
	
\begin{thm}[Bondy 1971] \label{Bondy} Let $G$ be a hamiltonian graph of order $p$.  If $G$ has at least $\dfrac{p^2}{4}$ edges, then $G$ is pancyclic  or $p$ is even and $G$ is the complete bipartite graph $K_{p/2,p/2}$.  
\end{thm}

In \cite{MM63}, Moon and Moser considered sufficient conditions for hamiltonicity in bipartite graphs. A bipartite graph $G$ of order $2n$ is {\it balanced}  if the vertex set of $G$ can be partitioned into two sets with $n$ vertices in each, such that every edge of $G$ joins vertices in different sets. If $G$ is a hamiltonian bipartite graph, then necessarily, $G$ is balanced.  The next theorem gives a sufficient condition for hamiltonicity in balanced bipartite graphs.

\begin{thm}[Moon and Moser 1963]\label{MoonMoser} Let $G$ be a balanced bipartite graph of order $p=2n$. If for every pair $u$ and $v$ of nonadjacent vertices in different partite sets, $d(u)+d(v)> n$, then $G$ is hamiltonian.
	
\end{thm}

Note that Theorem \ref{MoonMoser} improves the lower bound on degree sums in Theorem \ref{Ore} essentially from $p$ to $p/2$  if $G$ is a balanced bipartite graph.

\begin{cor}\label{corMM}
	Let $G$ be a balanced bipartite graph of order $p=2n \geq 4$. If $d(v) > n/2$ for every vertex $v$ of $G$, then $G$ is hamiltonian.
\end{cor}

No bipartite graph is pancyclic since bipartite graphs contain no odd cycles.  However, we can define a concept similar to pancyclicity for bipartite graphs.  We call a bipartite graph $G$ of order $2n$ {\it bipancyclic} if $G$ contains a $t$- cycle for every even integer $t$ between $4$ and $2n$. In \cite{ES88} Entringer and Schmeichel established an analogue to Theorem \ref{Bondy} for bipancyclicity. 

\begin{thm}[Entringer and Schmeichel 1988]  \label{ES} 
	Let $G$ be a balanced bipartite graph of order $p=2n \geq 4$.  If $G$ has at least $n^2 - n +2$ edges, then $G$ is bipancyclic.
\end{thm}
 
Quite recently, the result of Bondy in Theorem \ref{Bondy} was improved by Chen, Gould, Gu and Saito \cite{C18}. The improved result uses the concept of {\it chorded pancyclicity}, introduced by Cream, Gould and Hirohata \cite{C17}, which we recall now. \\

A {\it chord} of a cycle $C$ is an edge joining two non-consecutive vertices of $C$.  If a cycle $C$ of order $k$ has a chord, we call $C$  a {\it chorded k-cycle}.  A graph $G$ of order $p \ge 4$ is called {\it chorded pancyclic} if $G$ contains a chorded $k$-cycle for every integer $k$ with $4 \leq k \leq p$. As observed in \cite{C18} and \cite{C17}, chorded cycles are a fundamental tool for the study of the cycle distribution in a graph.  
 
 The following result by Chen et. al. appeared in \cite{C18} and will be used repeatedly later in this paper.

\begin{thm}[Chen, Gould, Gu and Saito 2018]\label{chorded18} Let $G$ be a hamiltonian graph of order $p$.  If $G$ has at least ${p^2}\over{4}$ edges, then $G$ is chorded pancyclic, or $p$ is even and $G=K_{p/2,p/2}$ or $G= K_3 \Box K_2$, the cartesian product of $K_3$ and $K_2$.  
\end{thm}

\smallskip

In Section 2, we present some definitions and known sufficient conditions for hamiltonicity in balanced $k$-partite graphs of order $kn$, for any integers $k\geq 3$ and $n\geq 1$. In Section 3 we prove that for all integers $k\geq 3$ and $n\geq 1$, every balanced $k$-partite graph  with $kn$ vertices and at least  ${{(k^2-k)n^2-2n(k-1)+4}\over 2}$ edges is chorded pancyclic, and in Section 4, we present a similar edge condition that guarantees chorded pancyclicity in $k$-partite graphs that are not necessarily balanced.

\bigskip

\section{Balanced $k$-partite graphs}

 A graph is {\it $k$-partite} if its vertex set can be partitioned into $k$ disjoint sets, or parts, in such a way that vertices in the same part are not adjacent. A $k$-partite graph is {\it balanced} if all its parts have the same number of vertices. A $k$-partite graph is {\it complete} if any two vertices in different parts are adjacent. The balanced complete $k$-partite graph of order $kn$, denoted $K_k(n)$ is the $k$-partite graph with $n$ vertices in each part, such that any two vertices in different parts are adjacent. Note that $K_k(1)$ is the complete graph of order $k$, also denoted by $K_k$. 
 
  Obviously, every graph $G$  can be viewed as a balanced $k$-partite graph of order $kn$ if we take $n=1$, and $k$  the order of $G$. The next theorem \cite{CJ97} and its corollary \cite{CFGJL95} extend Theorem \ref{MoonMoser} and Corollary \ref{corMM} to balanced $k$-partite graphs for $k \geq 3$.
  
  \bigskip
  
  \begin{thm}[Chen and Jacobson 1997] \label{degsumk} Let $k,n$ be integers, $k\geq 3$ and $n\geq 1$. Let $G$ be a balanced $k$-partite graph of order $p=kn$.

  	\textbf{Case 1.} If $k$ is even and  $d(u)+d(v) > \Big (k-\frac{4}{k+2} \Big ) n $ for every pair of nonadjacent vertices $u,v$ in different partite sets, then $G$ is hamiltonian.
  	
  	\smallskip
  	
  	\textbf{Case 2.}If $k$ is odd and   $d(u)+d(v) > \Big(k-\frac{2}{k+1}\Big ) n $ for every pair of nonadjacent vertices $u,v$ in different partite sets, then $G$ is hamiltonian.

  \end{thm}

\bigskip

\begin{cor}[Chen, Faudree, Gould, Jacobson and Lesniak 1995] \label{degreek} Let $k,n$ be integers, $k\geq 3$ and $n\geq 1$.  Let $G$ be a balanced $k$-partite graph of order $p=kn$.

		\textbf{Case 1.} If $k$ is even and  $d(v) > \Big ({k\over 2}-{2\over {k+2}} \Big ) n $ for every vertex $v$ of $G$, then $G$ is hamiltonian.
		
		\smallskip	
		
		\textbf{Case 2.} If $k$ is odd and  $d(v) > \Big ({k\over 2}-{1\over {k+1}} \Big ) n $  for every vertex $v$ of $G$, then $G$ is hamiltonian.		
		
	\end{cor}

\bigskip
	
Note Bondy's result in Theorem \ref{Bondy} cannot be applied in the conditions of Theorem \ref{degsumk} or Corollary \ref{degreek}. Hence, Theorem \ref{chorded18} cannot be applied either.
	
\section{Edge results for balanced $k$-partite graphs}

Corollary \ref{generalpancyclic} gives a sufficient edge condition for a graph to be pancyclic; Theorem \ref{chorded18} extends Corollary \ref{generalpancyclic} to a sufficient edge condition for a graph to be chorded pancyclic.  In \cite{A09}, Adamus gave a sufficient condition for a balanced tripartite graph to be pancyclic.

\bigskip
Since minimum degree at least $2$ is a necessary condition for a graph to be hamiltonian, Adamus noted that to guarantee that a balanced tripartite graph $G$ of order $3n$ is hamiltonian, we can remove at most $2n-2$ edges from the complete tripartite graph $K(n,n,n)$ to obtain $G$.  In other words, such a $G$ must have at least $3n^2-2n+2$ edges.  This condition is also sufficient.

\begin{thm}[Adamus 2009] \label{edgecondition3}
Let $G$ be a balanced tripartite graph of order $3n, n \geq 2$. If $G$ has at least $3n^2-2n+2$ edges, then $G$ is hamiltonian.
\end{thm}

As Adamus pointed out in \cite{A09}, while the edge condition in Corollary \ref{edgecondition1} follows directly from Ore's condition, the edge conditions in Theorem \ref{ES} and Theorem \ref{edgecondition3} follow from neither the Dirac minimum degree condition nor the Ore minimum degree sum condition. Adamus also noted that his edge condition for hamiltonicity does, in fact, give pancyclicity by Bondy's result in Theorem \ref{Bondy}. Hence, by Theorem \ref{chorded18}, the edge condition for hamiltonicity given by Adamus for balanced tripartite graphs actually gives chorded pancylicity.\\

In this section we give a sufficient edge condition for chorded pancyclicity in balanced $k$-partite  graphs of order $kn$ with  $k \geq 3$ and $n \geq 1$. Again, since minimum degree at least $2$ is necessary for hamiltonicity, we can remove at most $(k-1)n-2$ edges from the complete balanced $k$-partite graph $K_k(n)$ and still assure hamiltonicity.

\bigskip
The proof given by Adamus for $k=3$ relied only on Ore's sufficient condition (Theorem \ref{Ore}).  We include this case in our proof because the proof for all $k \geq 3$  follows rather quickly from the following classic theorem of P{\'o}sa \cite{P62}.  Furthermore, although Theorem ~\ref{ham} will follow from results in Section 4, we include its simple proof here. Understanding the proof of Theorem~\ref{ham} will help the reader follow the proof of Theorem ~\ref{hamnb}, which uses the same method but with additional nuances.
 
 \begin{thm} [P\'osa 1962]\label{Posa}
 	Let $G$ be a graph of order $p\geq 3$. If for every integer $r$, with $1\leq r <{p\over 2}$ the number of vertices of degree at most $r$ is less than $r$, then $G$ is hamiltonian.
 \end{thm}

\bigskip

\noindent {\bf Notation.} If $G=(V,E)$ is a graph and $S\subseteq V$, then $G[S]$ denotes the subgraph of $G$ induced by the vertices in $S$, i.e. $V( G[S])=S$ and $E(G[S])=\{ (x,y)\in E(G): x\in S, y\in S\}$. We use 
$||G||$ to denote the number of edges of $G$.\\

\begin{thm} \label{ham}
Let $k,n$ be integers, $k\geq 3$ and $n\geq 1$. Let $G$ be a balanced $k$-partite graph of order $p=kn$. If $G$ has at least $$ ||K_k(n)|| - ((k-1)n-2) = {{(k^2-k)n^2-2n(k-1)+4}\over 2} $$ edges, then $G$ is hamiltonian. 
\end{thm}

\bpf
We prove that $G$ satisfies P\'osa's condition by contradiction. If $G$ does not satisfy P\'osa's condition, there exists an integer $r$, $1\leq r <{{kn}\over 2}$, for which there are (at least) $r$ vertices $v_1,\ldots ,v_r$ such that $d_i=d_G(v_i)\leq r$.  Since, in fact, $G$ has minimum degree at least $2$, we can assume that $r \geq 2$. \smallskip

We can view $G$ as being obtained by deleting a set of edges from the complete $k$-partite graph $K_k(n)$.   Since in $K_k(n)$ every vertex has degree $(k-1)n$, to obtain vertices $v_1,\ldots ,v_r$ with degrees $d_1,\ldots ,d_r$ it is necessary to remove at least  $(k-1)n-d_i$  edges incident with vertex $v_i$, $i=1, \ldots ,r$. Then, the total number of removed edges is at least:
$$\sum_{i=1}^r \big ((k-1)n-d_i \big )-\big ( ||K_k(n)[\{v_1,\ldots ,v_r\}]||-||G[\{v_1,\ldots ,v_r\}]|| \big ) \geq \sum_{i=1}^r \big ((k-1)n-d_i \big ) -  ||K_k(n)[\{v_1,\ldots ,v_r\}]||$$
where the term $ ||K_k(n)[\{v_1,\ldots ,v_r\}]||-||G[\{v_1,\ldots ,v_r\}]||$ corresponds to the deleted edges that joined pairs of vertices in the set $\{v_1,\ldots ,v_r\}$ and were counted twice in the summation. \smallskip

The number $||K_k(n)[\{v_1,\ldots ,v_r\}]||$ depends on how the $r$ vertices are distributed among the $n$ parts. However, $$||K_k(n)[\{v_1,\ldots ,v_r\}]||\leq ||K_r||={{r(r-1)}\over 2},$$ so the number of edges removed from $K_k(n)$ to produce $G$ is at least: 
$$\sum_{i=1}^r \big ((k-1)n-d_i \big )- {{r(r-1)}\over 2}.$$ Since for every $i=1,\ldots ,r$, $d_i\leq r$, $$\sum_{i=1}^r \big ((k-1)n-d_i \big )- {{r(r-1)}\over 2}\geq  r \big ((k-1)n-r\big )-{{r(r-1)}\over 2}.$$ By assumption, at most $(k-1)n-2$ edges were removed from $K_k(n)$ to obtain $G$. It follows then, that

$$ r \big ((k-1)n-r\big)-{{r(r-1)}\over 2}\leq (k-1)n-2$$

%{\color{blue} $$ r (k-1)n-r^2-{{r(r-1)}\over 2}\leq (k-1)n-2$$}

\noindent  or equivalently, $$ r (k-1)n-r^2-{{r(r-1)}\over 2}< (k-1)n-1.$$

%{\color{blue} $$ (r-1) (k-1)n-r^2-{{r(r-1)}\over 2}< -1$$ $$ (r-1) (k-1)n-{{r(r-1)}\over 2}< r^2 -1$$ $$ (r-1) (k-1)n-{{r(r-1)}\over 2}< (r+1)(r-1)$$ $$ (r-1) \Big((k-1)n-{{r}\over 2}\Big)< (r+1)(r-1)$$}

\noindent Using some some basic arithmetic, this last inequality can be reduced to 

$$(r-1) \Big((k-1)n-{{r}\over 2}\Big) \mbox{  } <  \mbox{   } r^2 -1\mbox{  }=\mbox{  } (r+1)(r-1).$$

\noindent Since we are assuming $r\geq 2$, dividing both sides by $r-1$ we obtain 
$$(k-1)n-{{r}\over 2}< r+1$$ 

\noindent and this last inequality can be written as $(k-1)n-{{r}\over 2}\leq r$,  yielding  \begin{equation}\label{eqn1} 2(k-1)n\leq 3r. \end{equation}

\noindent Since $2\leq r <{{kn}\over 2}$, we have $$2(k-1)n \leq 3r< {{3kn}\over 2}$$ so,  $4(k-1)n <3kn$ and this implies $kn<4n$. This last inequality cannot hold if $k\geq 4$. Therefore, if $k\geq 4$ then $G$ is hamiltonian by P\'osa's condition.

\smallskip

In the case $k=3$ (Adamus' result), we first show that P\'osa's condition holds for $r=2,3$. From (1) we know $ 2(k-1)n\leq 3r$ and since $k=3$ in we obtain $4n\leq 3r$. Besides, since $r<{{3n}\over 2}$, it must be $4n\leq 3r<{{9n}\over 2}$ and the two inequalities are not compatible. If $r=2$, the leftmost equality is only possible for $n=1$ but the rightmost inequality only holds for $n\geq 2$; if $r=3$ the leftmost inequality holds if $n=1,2$ but the rightmost inequality only holds for $n>3$.

%{\color{blue} In this case $r=2:$ we have $4n\leq 6<{{9n}\over 2}$ or $8n\leq 12<9n$, but $8n\leq 12$ means $n<{3\over 2}$ or $n\leq 1$ and $12<9n$ means ${4\over 3}<n$ or $n\geq 2$. }

%{\color{blue} If $r=3$ we obtain $4n\leq 9<{{9n}\over 2}$ or $8n\leq 18<9n$. But $8n\leq 18$ implies  $n=1$ or $n=2$, while $18<9n$ implies $n>3$. }

Let us now prove P\'osa's condition for $r\geq 4$.  By contradiction, assume there exists $r$, $4\leq r<{{3n}\over 2}$ such that $G$ has $r$ vertices $v_1,\ldots ,v_r$ such that $d_i=d_G(v_i)\leq r$. Then, $r\geq 4$ and the fact that there can be at most $5$ edges between any four vertices in $K_3(n)$, together imply that the number of edges deleted from $K_3(n)$ to create $G$ is at least, 
$$2n-d_1+2n-d_2+2n-d_3+2n-d_4-5\geq 8n-4r-5$$
%{ \color{blue} Since $d_i=d_G(v_i)\leq r$ for $i=1,2,3,4$, then $$2n-d_1+2n-d_2+2n-d_3+2n-d_4-5 \geq 2n-r+2n-r+2n-r+2n-r-5= 8n-4r-5.$$ Imposing $2n-2\geq 2n-d_1+2n-d_2+2n-d_3+2n-d_4-5$, $$2n-2\geq 2n-d_1+2n-d_2+2n-d_3+2n-d_4-5\geq 2n-r+2n-r+2n-r+2n-r-5= 8n-4r-5$$}\\
\noindent and the condition on the size of $G$ guarantees that $8n-4r-5\leq 2n-2$,  which implies ${{3n}\over 2}\leq r$ and contradicts the condition 
 $ r <{{kn}\over 2}$. %{\color{blue} From $8n-4r-5\leq 2n-2$ we obtain $6n\leq 4r+3$ or $6n<4r+4=4(r+1)$ so ${{6n}\over 4}<r+1$ or ${{3n}\over 2}\leq r <{{kn}\over 2}={{3n}\over 2}$}
 \epf

The following result is a direct consequence of Theorem \ref{ham} combined with Theorem \ref{chorded18}.

\begin{cor}\label{pan}
Let $k,n$ be integers, $k\geq 3$ and $n\geq 1$. Let $G$ be a balanced $k$-partite graph of order $p=kn$. If $G$ has at least $ {{(k^2-k)n^2-2n(k-1)+4}\over 2}$ edges, then $G$ is chorded pancyclic. 
\end{cor}

\bpf
 Observe that $G$ is neither $K_{{p/2},{p/2}}$ nor $K_3 \Box K_2$. Then, by Theorem \ref{chorded18}, since $G$ has order $p=kn$, it is sufficient to show that
$${{(k^2-k)n^2-2n(k-1)+4}\over 2}\geq {{(kn)^2}\over 4}$$ or equivalently, $${{(k^2-k)n^2-2n(k-1)+4}\over 2}> {{(kn)^2}\over 4}-1.$$ Since ${{(kn)^2}\over 4}-1= {{(kn-2)(kn+2)}\over 4}$, this inequality can be written as
\begin{equation}\label{eq:2}{(k^2-k)n^2-2n(k-1)+4}> {{(kn-2)(kn+2)}\over 2}.\end{equation} 
It is straightforward to verify that  \ref{eq:2} holds if $k=4$, $n=1$ and if $k=3$, $n=1,2$. Assume that this is not the case. We then show 
\begin{equation}\label{eq:3}(k^2-k)n^2-2n(k-1)>{{(kn-2)(kn+2)}\over 2}\end{equation} 
which suffices to complete the proof.  Using basic arithmetic it can be shown that $$(k^2-k)n^2-2n(k-1)={(k-1)n(nk-2)}$$ 
%{\color{blue} $$(k^2-k)n^2-2n(k-1)=k(k-1)n^2-2n(k-1)=n(k-1)(nk-2)$$ Then, $$ (k^2-k)n^2-2n(k-1)=n(k-1)(nk-2)>{{(kn-2)(kn+2)}\over 2}$$ and since $kn-2\not= 2$ (because $k\geq 3$, $n\geq 1$) we divide by $kn-2$ and obtain $n(k-1)>{{(kn+2)}\over 2}$ }
so the previous inequality is equivalent to 
$${2(k-1)n}> kn+2$$ %{\color{blue} $$2kn-2n>kn+2$$ $$kn=2kn-kn >2n+2=2(n+1)$$}
and can be reduced to $kn>2(n+1)$. %{\color{blue} Since $k\geq 3$, then  $kn\geq 3n =2n+n>2n+2$ if $n>2$. If $n=2$, we have $2k>6$ and the only case when this does not happen is if $k=3$. If $n=1$ we have $k>4$ and this fails for $k=3,4$} 
This inequality holds for any $k\geq 3$ and $n\geq 1$, except when $k=4$ and $n=1$ or  $k=3$ and $n=1,2$. This completes the proof.

\epf
\bigskip

The previous two results are best possible since the graph obtained from $K_k(n)$ by removing all but one edge from any vertex gives a nonhamiltonian graph with exactly ${{(k^2-k)n^2-2n(k-1)+4}\over 2} - 1$ edges.

\bigskip

Analogosuly to the edge conditions for bipartite graphs in Theorem \ref{ES} and for tripartite graphs in Theorem \ref{edgecondition3}, the edge condition in Theorem \ref{ham} follows neither from the Dirac minimum degree condition nor the Ore minimum degree sum condition.

\bigskip
We close by noting that the number of edges required for bipancyclicity in Theorem 10 is that of Theorem 14 if $k$ is replaced with $2$.

\section{Edge results for general $k$-partite graphs}

We begin by setting up the notation needed to study general $k$-partite graphs.

\bigskip

\noi {\bf Notation.} For an integer $k\geq 3$, consider a $k$-tuple of positive integers $(n_1,\ldots ,n_k)$ such that $n_1\geq n_2\geq \ldots \geq n_k$. Define $p=\sum _{i=1}^k{n_i}$ and let ${\cal G}(n_1,\ldots ,n_k)$ denote the set of all $k$-partite graphs with parts $V_1,\ldots ,V_k$ such that $|V_i|=n_i$ for every $i=1,\ldots ,k$. Note that, as in the previous sections, $p$ denotes the order of $G$.

\bigskip
\smallskip

If $n_1> {p\over 2}$, then no graph in ${\cal G}(n_1,\ldots ,n_k)$ is hamiltonian, so a necessary condition for our work is that $n_1\leq {p\over 2}$. The graphs in ${\cal G}(n_1,\ldots ,n_k)$ are the result of removing edges from the complete $k$-partite graph $K(n_1,\ldots ,n_k)$. The condition $n_1\leq {n\over 2}$ guarantees that $K(n_1,\ldots ,n_k)$ is hamiltonian, and we want to determine the maximum integer $m$ such that removing any set of at most $m$ edges from  $K(n_1,\ldots ,n_k)$ yields a hamiltonian graph. 

\smallskip

Another necessary condition for a graph to have a hamiltonian cycle is that every vertex must have at least degree $2$. In the graph $K(n_1,\ldots ,n_k)$, each vertex in $V_i$ has degree $p-n_i$, for $i=1,\ldots ,k$. Thus, the condition $n_1\geq n_2\geq \ldots \geq n_k$ implies that the minimum degree of $K(n_1,\ldots ,n_k)$ is $p-n_1$. Therefore, a necessary condition for the integer $m$ that we want to determine, is that $m\leq p-n_1-2$.

\bigskip

The following results show that if $n_1\leq {p\over 2}$, any graph obtained by deleting at most $p-n_1-2$ edges from $K(n_1,\ldots ,n_k)$ is hamiltonian. As a consequence, these two necessary conditions for hamiltonicity turn out to be sufficient. Our first theorem corresponds to the case when $n_1< {p\over 2}-1$ and its proof follows from P\'osa's condition for hamiltonicity as in the balanced case. 

\bigskip

\begin{thm}  \label{hamnb}
Let $k\geq 3$ be an integer and let $(n_1,\ldots ,n_k)$ be a $k$-tuple of positive integers such that $n_1\geq n_2\geq \ldots \geq n_k$. Let $p=\sum _{i=1}^k{n_i}$. If $n_1< {p\over 2}-1$, then every graph $G$ in ${\cal G}(n_1,\ldots ,n_k)$ with at least 

$$ ||K(n_1,\ldots, n_k)|| - (p-n_1-2) = {{4 -2p+2n_1+\sum _{i=1}^{k}  n_i(p-n_i)}\over 2}$$ edges is hamiltonian. 
\end{thm}

\bpf
We prove that $G$ satisfies P\'osa's condition for hamiltonicity by contradiction, as we did in the proof of Theorem \ref{ham}. If $G$ does not satisfy P\'osa's condition, there exists an integer $r$, $1\leq r <{p\over 2}$, for which there are (at least) $r$ vertices $v_1,\ldots ,v_r$ such that $d_G(v_j)\leq r$.  As in the proof of Theorem \ref{ham}, we may assume $r \geq 2.$ 

\smallskip
Furthermore, if there exists a vertex $v$ with $d_G(v)=2$ it is necessary for $v$ to have minimum degree in $K(n_1,\ldots , n_k)$, and also that each of the edges removed from $K(n_1, \ldots , n_k)$ to produce $G$ is incident with $v$. The only way to obtain a second vertex $u$ with $d_G(u)=2$ is if there exists a neighbor of $v$ with degree $3$ in $K(n_1, \ldots , n_k)$. However, this can only happen if $p=4$, but in this case $r<2$. Thus, we may assume $r\geq 3$.

As in the proof of Theorem \ref{ham} every graph $G$ in ${\cal G}(n_1,\ldots ,n_k)$ is obtained by deleting some edges from $K(n_1,\ldots ,n_k)$. If a vertex $v$ has $d_G(v)\leq r$, at least $p-n_1-r$ edges incident with $v$ were removed from $K(n_1,\ldots ,n_k)$. Since at most  $p-n_1-2$ edges were removed from $K(n_1,\dots ,n_k)$ to produce $G$, it must be: $$r \big (p-n_1-r \big )- {{r(r-1)}\over 2} < p-n_1-1,$$ 
%{\color{blue} $$r \big (p-n_1\big )-r^2- {{r(r-1)}\over 2} < (p-n_1)-1$$ $$(r-1) \big (p-n_1\big )-{{r(r-1)}\over 2} < r^2 -1 $$}
or equivalently, 
$$(r-1)(p-n_1)-{{r(r-1)}\over 2} \mbox{  } < \mbox{  } r^2-1 \mbox{  } = \mbox{  } (r+1)(r-1).$$ Then, $(r-1)(p-n_1)-{{r(r-1)}\over 2} < (r+1)(r-1),$ and  since $r\geq 3$, dividing by $r-1$, we obtain $$p-n_1-{{r}\over 2} < r+1.$$ %{\color{blue} $$p-n_1-{{r}\over 2} \leq r$$ $$p-n_1\leq{{3r}\over 2} $$}
As in balance case, this inequality can be reduced to  \begin{equation}\label{eq:4}
2(p-n_1)\leq  {{3r}} \end{equation} where $p-n_1$ in \ref{eq:4} is the same as $(k-1)n$ in \ref{eqn1}.

 Since $r<{p\over 2}$, from equation \ref{eq:4} we conclude $2(p-n_1)< {{3p}\over 2}$ and using some basic arithmetic this expression can be reduced to $n_1>{p\over 4}$. Therefore, when  $n_1\leq {p\over 4}$  P\'osa's condition guarantees that $G$ is hamiltonian.\smallskip

In the case ${p\over 4}<n_1< {p\over 2}-1$, let $E'$ be the set of all $t$ edges removed from $K(n_1,\dots ,n_k)$ to produce $G$, and assume
$$E'=\{a_ib_i: 1\leq i\leq t\} \mbox{ with } t\leq p-n_1-2.$$ Define $$V'=\{v\in V(K(n_1,\dots ,n_k)): \exists i, 1\leq i\leq t, \mbox{ such that } v=a_i \mbox{ or } v=b_i \},$$ so that  $G'=(V',E')$ is the subgraph of $K(n_1,\dots ,n_k)$ induced by the edges removed from $K(n_1,\dots ,n_k)$ to produce $G$. 

\smallskip

If $G$ contains $r$ vertices of degree at most $r$, then $G'$ must contain at least $r$ vertices of degree at least $p-n_1-r$. Thus \begin{equation}\label{eq:5} \sum _{u\in V'} d_{G'}(u)\geq r(p-n_1-r). \end{equation}
 
 \smallskip

\noindent At the same time, $\sum _{u\in V'} d_{G'}(u)=2||G'||$ and $||G'||=|E'|\leq p-n_1-2$, so it must be \begin{equation}\label{eq:6}
\sum _{u\in V'} d_{G'}(u) \leq 2(p-n_1-2). \end{equation}

\noindent  Combining equations \ref{eq:5} and \ref{eq:6} we obtain $$2(p-n_1-2)\geq \sum _{u\in V'} d_{G'}(u)\geq r(p-n_1-r).$$

\noindent As a consequence, $2(p-n_1-2)\geq r(p-n_1-r)$, and this expression can be rewritten as %{\color{blue} $$2(p-n_1)-4\geq r(p-n_1)-r^2$$  $$r^2-4\geq (r-2)(p-n_1)$$  } 
$$(r-2)(r+2) \mbox{  } = \mbox{  } r^2-4 \mbox{  } \geq \mbox{  } (r-2)(p-n_1)$$

Using that $r\geq 3$, this expression can be reduced to \begin{equation}\label{eq:7} r+2\geq p-n_1,\end{equation} and adding the condition $r<{p\over 2}$,  $${p\over 2}+2> r+2\geq p-n_1.$$ % {\color{blue} $${p\over 2}+2>  p-n_1 $$  $${p\over 2}+1\geq  p-n_1 $$   $$n_1\geq  p-{p\over 2}-1={p\over 2}-1$$}

 From this expression we obtain $n_1\geq  {p\over 2}-1,$ which contradicts $n_1< {p\over 2}-1$. \epf
 
 \bigskip
 We now consider the case ${p\over 2}-1\leq n_1\leq {p\over 2}$. Depending on $p$ being even or odd, there are three cases where this can happen:
 
 \begin{itemize}
 \item [1)] $p$ even and $n_1={p\over 2}$,
 \item [2)] $p$ odd and $n_1={ {p-1}\over 2}$ and
 \item [3)] $p$ even and $n_1={p\over 2}-1.$
 \end{itemize}
 
 \begin{rem}
 
The technique we used to prove Theorem \ref{ham} and Theorem \ref{hamnb} cannot be applied if ${p\over 2}-1\leq n_1\leq {p\over 2}$. Indeed, the following examples show that for each of the cases above, it is possible to construct at least one family of graphs satisfying the edge the condition in Theorem \ref{hamnb} for which P\'osa's condition does not hold. However, the graphs in the families we present are hamiltonian.  

 \begin{itemize}

\item[1)]  Assume $p$ is even and $n_1={p\over 2}$.

For any integer $a\geq 5$, let $G$ be the graph in ${\cal G}(a,a-2,2)$ obtained by choosing any $a-2$ vertices $v_1,\ldots ,v_{a-2}$ among the $a$ vertices in $V_1$, a vertex $u$ in $V_2$, and removing the $a-2$ edges $v_iu$, for every $i=1,\ldots ,a-2$. Then, $G$ satisfies the edge condition in Theorem \ref{hamnb} but $G$ fails P\'osa's condition for $r=a-1<{p\over 2}$.
 
\item[2)]  Assume $p$ is odd and $n_1={ {p-1}\over 2}$.

For any integer $a\geq 3$, let $G$ be the graph in ${\cal G}(a,a-1,2)$ obtained by choosing $a-1$ vertices $v_1,\ldots ,v_{a-1}$ among the $a$ vertices in $V_1$, a vertex $u$ in $V_2$, and removing the $a-1$ edges $v_iu$, for $i=1,\ldots ,a-1$. Then,  $G$ satisfies the edge condition in Theorem \ref{hamnb} but $G$ fails P\'osa's condition for $r=a<{p\over 2}$.

\item[3)]  Assume $p$ is even and $n_1={p\over 2}-1$.

For an integer $a\geq 1$, consider $K(4a,4a,2)$ and assume $V_1=\{u_1^1,\ldots ,u_a^1\}\cup \{u_1^2,\ldots ,u_a^2\}\cup \{u_1^3,\ldots ,u_a^3\}\cup \{u_1^4,\ldots ,u_a^4\}$ and $V_2=\{v_1^1,\ldots ,v_a^1\}\cup \{v_1^2,\ldots ,v_a^2\}\cup \{v_1^3,\ldots ,v_a^3\}\cup \{v_1^4,\ldots ,v_a^4\}$.  Let $G$ be the graph in ${\cal G}(4a,4a,2)$ obtained by removing  from $K(4a,4a,2)$ the $4a$ edges $u_i^1v_i^1$,  $u_i^2v_i^2$,  $u_i^1v_i^2$ and  $u_i^2v_i^1$ for $i=1,\ldots ,a$. Then,  $G$ satisfies the edge condition in Theorem \ref{hamnb} but $G$ fails P\'osa's condition for $r=4a={p\over 2}-1$.

\end{itemize}
\end{rem}

Next, we prove the edge condition when ${p\over 2}-1\leq n_1\leq {p\over 2}$. In the cases when  P\'osa's condition is not satisfied, we apply Theorem \ref{ES} to a balanced complete bipartite subgraph of $K(n_1, \ldots , n_k)$. \\
 
In the case $n$ even and $n_1= {p\over 2}$, since exactly half of the vertices are in $V_1$, if there is a hamiltonian cycle in a graph $G$ in ${\cal G}(n_1,\ldots ,n_k)$, then every edge in the cycle must have an endpoint in $V_1$ and the other in $V\setminus V_1$. Therefore, edges having both endpoints in $V\setminus V_1$ do not affect the hamiltonicity of $G$ and we can prove a stronger result in this case.
 
\begin{thm} \label{evenhalfthm}
Let $k\geq 3$ be an integer and let $(n_1,\ldots ,n_k)$ be a $k$-tuple of positive integers $n_1\geq n_2\geq \ldots \geq n_k$. Let $p=\sum _{i=1}^k{n_i}$. If $p$ is even and $n_1= {p\over 2}$, then every graph $G$ in ${\cal G}(n_1,\ldots ,n_k)$ having at least $\big( {p\over 2}\big)^2- {p\over 2}+2$ edges with an endpoint in $V_1$, is hamiltonian.
\end{thm}

\bpf
Assume that all edges with both endpoints in $V\setminus V_1$ are deleted from $K(n_1, \ldots , n_k)$ and as a result, $G$ is obtained by deleting at most $p-n_1-2$ edges from $K_2\big ({p\over 2}\big )$. Then, $G$ has at least $||K_2\big ({p\over 2}\big )||-(p-n_1-2)=\big( {p\over 2}\big)^2- {p\over 2}+2$ edges, and by Theorem \ref{ES} $G$ is hamiltonian.
\epf

\begin{cor} \label{evenhalf}
Let $k\geq 3$ be an integer and let $(n_1,\ldots ,n_k)$ be a $k$-tuple of positive integers $n_1\geq n_2\geq \ldots \geq n_k$. Let $p=\sum _{i=1}^k{n_i}$. If $p$ is even and $n_1= {p\over 2}$, then every graph $G$ in ${\cal G}(n_1,\ldots ,n_k)$ with at least 

$$ ||K(n_1, \ldots , n_k)|| - (p-n_1-2) = { {4 -2p+2n_1+\sum _{i=1}^{k}  n_i(p-n_i)}\over 2}$$

edges is hamiltonian. 
\end{cor}

\bpf
Since the graph $G$ is in ${\cal G}(n_1,\ldots ,n_k)$ and $||G||\geq ||K(n_1, \ldots , n_k)|| - (p-n_1-2)$, then $G$ is the result of deleting at most $p-n_1-2={p\over 2}-2$ edges from $K(n_1, \ldots , n_k)$. As a consequence, even if all the ${p\over 2}-2$ deleted edges are selected from the $n_1(n-n_1)=\big( {p\over 2}\big)^2$ edges  incident with a vertex in $V_1$, $G$ has at least $\big( {p\over 2}\big)^2-{p\over 2}+2$ edges with an endpoint in $V_1$, so Theorem \ref{evenhalfthm} guarantees that $G$ is hamiltonian.
\epf

 \bigskip

In the case when $p$ is odd and $n_1={{p-1}\over 2}$, we have $|V\setminus V_1|=|V_1|+1$. Therefore, if there is a hamiltonian cycle in in a graph $G$ in ${\cal G}(n_1,\ldots ,n_k)$, then there is exactly one edge in the cycle having both endpoints in $V\setminus V_1$. 

\begin{thm} \label{oddhalf}
Let $k\geq 3$ be an integer and let $(n_1,\ldots ,n_k)$ be a $k$-tuple of positive integers $n_1\geq n_2\geq \ldots \geq n_k$. Let $p=\sum _{i=1}^k{n_i}$. If $p$ is odd and $n_1= {{p-1}\over 2}$, then every graph $G$ in ${\cal G}(n_1,\ldots ,n_k)$ with at least 

$$ ||K(n_1, \ldots ,n_k)|| - (p-n_1-2) = { {4 -2p+2n_1+\sum _{i=1}^{k}  n_i(p-n_i)}\over 2}$$

edges is hamiltonian. 
\end{thm}

\smallskip

\bpf
Since $G$ has at least $||K(n_1,  \ldots ,n_k)|| - (p-n_1-2)$ edges, $G$ is the result of deleting at most $p-n_1-2={{p-3}\over 2}$ edges from $K(n_1, \ldots , n_k)$. 

In $K(n_1, \ldots ,n_k)$ there are ${{p+1}\over 2}$ vertices in  $V\setminus V_1$. Therefore, in every graph $G$ obtained by removing at most ${{p-3}\over 2}$ edges from $K(n_1, \ldots , n_k)$,  there exists a vertex  $w$ in $V\setminus V_1$ such that $V_1\subseteq N_G(w)$. Define $G'=G-w$, and observe that since $G$ is obtained by removing at most ${{p-3}\over 2}$ from $K(n_1, \ldots , n_k)$, then $G'$ is the result of removing at most ${{p-3}\over 2}$ edges from $K(n_1, \ldots , n_k)-w$.  \\

Let us distinguish two types of edges in $K(n_1, \ldots , n_k)-w$. Edges of {\it type 1} are those with an endpoint in $V_1$ and the other in $V\setminus V_1$, while edges of {\it type 2} are those with both endpoints in $V\setminus V_1$. \\

First, consider the case when $G'$ is a graph obtained by removing from $K(n_1, n_2, \ldots ,n_k)-w$ at most ${{p-5}\over 2}$ edges of type $1$ and at most ${{p-3}\over 2}$ edges of type $2$.  \smallskip

Since $k\geq 3$, the graph $K(n_1, n_2, \ldots ,n_k)-w$ has at least ${{p-1}\over 2}$ edges of type $2$. Thus, even when deleting exactly ${{p-3}\over 2}$ edges of type $2$,  there exists an edge $e$ in $G'$ with both endpoints in $V\setminus V_1$.   \smallskip

Observe that $K_2\big({{p-1}\over 2}\big)$, the balanced complete bipartite graph of order $p-1$, is a spanning subgraph of $K(n_1, \ldots , n_k)-w$ and it contains all edges of type $1$ in $K(n_1, \ldots , n_k)-w$. As a consequence, there is a spanning subgraph of $G'$ that results from deleting at most ${{p-5}\over 2}$ edges from $K_2\big({{p-1}\over 2}\big)$. By Theorem \ref{ES}, removing at most ${{p-5}\over 2}$ edges from $K_2\big({{p-1}\over 2}\big)$ yields a hamiltonian graph. Thus, $G'$ has a hamiltonian spanning subgraph so $G'$ is hamiltonian.  \smallskip

We show next that if $G'$ is hamiltonian, then $G$ is also hamiltonian. To do this, we construct a hamiltonian cycle in $G$ from a hamiltonian cycle in $G'$, together with $w$ and the edge $e$ from above. Observe that $G'$ is a $k$-partite graph of even order $p-1$ with exactly ${{p-1}\over 2}$ vertices in $V_1$. Therefore, a hamiltonian cycle in $G'$ is an alternating sequence of vertices in $V_1$ and vertices in $V\setminus V_1$. Let $C=v_1,\ldots ,v_{p-1},v_1$ be a hamiltonian cycle in $G'$, where the vertices with odd sub-indices are in $V_1$ and the vertices with even sub-indices are in $V\setminus V_1$. Since $e$ has both endpoints in $V\setminus V_1$, there exist integers $i$ and $j$, $1\leq i<j\leq p-1$ such that $e=v_{2i}v_{2j}$. Then, $v_{2i-1}\in V_1$ and $v_{2j-1}\in V_1$  and $V_1\subseteq N_G(w)$ implies that there is a path $v_{2i-1},w,v_{2j-1}$ in $G$, so $v_1,\ldots ,v_{2i-1},w,v_{2j-1},v_{2j-2}, \ldots , v_{2i}, v_{2j}, \ldots ,v_{p-1}$ is a hamiltonian cycle in $G$. \\

Now consider the case when $G'$ is obtained by deleting from $K(n_1,\ldots ,n_k)-w$, exactly ${{p-3}\over 2}$ edges of type $1$. In this case, since no edges of type $2$ are removed, we have $N_G(w)=N_{K(n_1,\ldots ,n_k)}(w)$.  \smallskip

As in the previous case, $K_2\big({{p-1}\over 2}\big)$ is a spanning subgraph of $K(n_1, \ldots , n_k)-w$ and it contains all edges  of type $1$ in $K(n_1, \ldots , n_k)-w$. Then, there is a spanning subgraph of $G'$ that results from deleting at most ${{p-3}\over 2}$ edges from $K_2\big({{p-1}\over 2}\big)$. By Theorem \ref{ES}, removing at most ${{p-5}\over 2}$ edges from $K_2\big({{p-1}\over 2}\big)$ yields a hamiltonian graph, so we conclude that $G'$ has a hamiltonian path $P$.  Assume  $P=v_1,\ldots ,v_{p-1}$ where the vertices with odd sub-indices are in $V_1$ and the vertices with even sub-indices are in $V\setminus V_1$. If  $v_{p-1}$ and $w$ are different parts, then the path $v_1,w,v_{p-1}$ together with $P$ form a hamiltonian cycle in $G$. If $v_{p-1}$ and $w$ are in the same part,  since $k\geq 3$ and no edges of type $2$ had been removed,  $v_{p-1}$ has a neighbor in $V\setminus V_1$, say $v_{2i}$. Since $P$ alternate vertices in  $V\setminus V_1$ and vertices in $V_1$, $v_{2i+1}$ is in $V_1\subseteq N_G(w)$. Then, $v_{p-1},v_{2i},v_{2i-1}, \ldots ,v_1,w,v_{2i+1},\ldots ,v_{p-1}$ is a hamiltonian cycle in $G$. 
 \epf

  \bigskip
  
 \begin{thm} \label{evenhalf-1}
Let $k\geq 3$ be an integer and let $(n_1,\ldots ,n_k)$ be a $k$-tuple of positive integers $n_1\geq n_2\geq \ldots \geq n_k$. Let $p=\sum _{i=1}^k{n_i}$. If $p$ is even and $n_1= {p\over 2}-1$, then every graph $G$ in ${\cal G}(n_1,\ldots ,n_k)$ with at least 

$$ ||K(n_1,\ldots ,n_k)|| - (p-n_1-2) = { {4 -2p+2n_1+\sum _{i=1}^{k}  n_i(p-n_i)}\over 2}$$

  edges is hamiltonian. 
\end{thm}

\bpf
If $n_2<n_1$ we will show that $G$ satisfies P\'osa's condition. By contradiction, assume there exists an integer $r$, $1\leq r< {p\over 2}$ for which there exist $r$ vertices $v_1,\ldots ,v_r$ with $d_G(v_i)\leq r$.  As in the proof of Theorem \ref{hamnb}, this implies $r+2\geq p-n_1$ \ref{eq:7}. Replacing $n_1= {p\over 2}-1$ we obtain %$r+2\geq {n\over 2}+1$ so 
$r\geq {p\over 2}-1$ and thus, P\'osa's condition holds for $r\leq {p\over 2}-2$.  Since $r< {p\over 2}$, the only remaining possibility is $r={p\over 2}-1$. However, if $r={p\over 2}-1$, then  %$n-n_1-r=n-{n\over 2}+1-{n\over 2}+1=2$. 
$p-n_1-r=2$ and for each vertex $v$ in $V_1$ with $d_G(v)\leq r$, it is necessary to delete from $K(n_1,\ldots ,n_k)$ at least two edges incident with $v$. The assumption $n_2<n_1$  implies $n_2\leq {p\over 2}-2$, so $p-n_2-r\geq 3$ and for each vertex $u$ in $V\setminus V_1$ with $d_G(u)\leq r$, it is necessary to remove at least 
three edges incident with $u$ in $K(n_1,\ldots ,n_k)$. As a result, there are at most $\big\lfloor  {{p-2} \over 4}\big\rfloor + \big\lfloor  {{p-2} \over 6}\big\rfloor< r$ vertices of degree at most $r$ in $G$ so P\'osa's condition also holds for $r={p\over 2}-1$ and $G$ is hamiltonian. \\

If $n_2=n_1$, then $n_1+n_2=p-2$. Thus, $G$ is obtained by removing at most $p-n_1-2= {p\over 2}-1$ edges from $K( {p\over 2}-1, {p\over 2}-1,2)$ or $K\big({p\over 2}-1, {p\over 2}-1,1,1\big)$. Since $K\big({p\over 2}-1, {p\over 2}-1,1,1\big)$ has one more edge than $K( {p\over 2}-1, {p\over 2}-1,2)$, it is sufficient to show that any graph $G$ obtained by removing ${p\over 2}-1$ edges from $K\big( {p\over 2}-1, {p\over 2}-1,2\big)$ is hamiltonian.  \smallskip

Note that $G[V_1\cup V_2]$ is a sub-graph of $K_2\big( {p\over 2}-1\big)$.  Then, Theorem \ref{ES} guarantees that if at most ${p\over 2}-3$ of the edges removed from $K\big( {p\over 2}-1, {p\over 2}-1,2\big)$ join a vertex in $V_1$ and a vertex in $V_2$, then $G[V_1\cup V_2]$ is hamiltonian. Therefore, there exists a cycle in $G$ that contains all $p-2$ vertices in $V_1\cup V_2$. Denote such cycle as  $C=v_1,\ldots ,v_{p-2},v_1$, where the vertices with odd sub-indices are in $V_1$, the vertices with even sub-indices are in $V_2$ and assume $V_3=\{x,y\}$. If there exist two different edges $v_{2i}v_{2i+1}$ and $v_{2j}v_{2j+1}$ in $C$ such that $v_{2i}$, $v_{2i+1}$ are in $N_G(x)$ and $v_{2j}$, $v_{2j+1}$ are in $N_G(y)$, then respectively replacing these edges in $C$ with $v_{2i},x,v_{2i+1}$ and $v_{2j},y,v_{2j+1}$ we obtain a hamiltonian cycle in $G$. \smallskip

Since $G$ results from deleting ${p\over 2}-1$ edges from $K\big( {p\over 2}-1, {p\over 2}-1,2\big)$,  $d_G(x)+d_G(y)\geq 3 \big({p\over 2}-1 \big)$. If $d_G(x)={p\over 2}-1$, then $d_G(y)=p-2$. If $x$ has neighbors in $V_1$ and $V_2$, then there exists an edge $v_{2i}v_{2i+1}$ with  $v_{2i}$, $v_{2i+1}$ in $N_G(x)$, and for any other edge $v_{2j}v_{2j+1}$  we have $v_{2j}$,$v_{2j+1}$ in $N_G(y)$, so we construct a hamiltonian cycle in $G$ as above. If $N_G(x)=V_1$, then $v_1,x,v_3,v_2,y,v_4, \ldots ,v_{p-2},v_1$ is a hamiltonian cycle in $G$; if $N_G(x)=V_2$ we proceed in the same way as when $N_G(x)=V_1$. If $d_G(y)={p\over 2}-1$, then $d_G(x)=p-2$, so we proceed as in the previous case. In all other cases, there exist two different edges $v_{2i}v_{2i+1}$ and $v_{2j}v_{2j+1}$ in $C$ such that $v_{2i}$,$v_{2i+1}$ are in $N_G(x)$ and $v_{2j}$, $v_{2j+1}$ are in $N_G(y)$, so we construct a hamiltonian cycle in $G$. \smallskip

Now, consider that ${p\over 2}-2$ of the edges removed from $K\big( {p\over 2}-1, {p\over 2}-1,2\big)$, join a vertex in $V_1$ and a vertex in $V_2$. After removing the first ${p\over 2}-3$ edges, there is a cycle $C=v_1,\ldots ,v_{p-2},v_1$ containing the $p-2$ vertices in $V_1\cup V_2$. Suppose that an edge $v_{2i}v_{2i+1}$ in $C$ is deleted. Since at most one additional edge can be removed,  at least one  of the two vertices in $V_3$ remain adjacent with all vertices in $V_1\cup V_2$. Suppose $N_G(x)=V_1\cup V_2$. Then, $v_1,\ldots ,v_{2i},x,v_{2i+1},\ldots ,v_{p-2},v_1$ is a hamiltonian cycle in $G$. \smallskip

Finally, if two edges in the hamiltonian cycle are removed, say $v_{2i}v_{2i+1}$ and $v_{2j}v_{2j+1}$ with $i<j$, then $N_G(x)=V_1\cup V_2$, $N_G(y)=V_1\cup V_2$, and $v_1,\ldots ,v_{2i},x,v_{2i+1},\ldots ,v_{2j},y,v_{2j+1},\ldots ,v_{p-2},v_1$ is a hamiltonian cycle in $G$. 
 \epf
\smallskip
 Since every balanced $k$-partite graph satisfies the conditions of one of Theorems \ref{hamnb}, \ref{evenhalf}, \ref{oddhalf} or \ref{evenhalf-1}, Theorem \ref{ham} is a corollary of the results in this section.  
  \bigskip
  
 Analogously to the balanced case, combining the hamiltonicity results in Theorems \ref{hamnb}, \ref{evenhalf}, \ref{oddhalf} and \ref{evenhalf-1} with the minimum edge condition for chorded pancyclicity in Theorem \ref{chorded18} we obtain the primary result of this section.

\begin{thm} 
Let $k\geq 3$ be an integer and let $(n_1,\ldots ,n_k)$ be a $k$-tuple of positive integers $n_1\geq n_2\geq \ldots \geq n_k$. Let $p=\sum _{i=1}^k{n_i}$. If  $n_1\leq  {{p}\over 2}$, then every graph $G$ in ${\cal G}(n_1,\ldots ,n_k)$ having at least 

$$ ||K(n_1,  \ldots ,n_k)|| - (p-n_1-2) = { {4 -2p+2n_1+\sum _{i=1}^{k}  n_i(p-n_i)}\over 2}$$

edges is chorded pancyclic. 
\end{thm}

\bpf
By Theorem \ref{chorded18}, it is sufficient to show $||K(n_1,  \ldots ,n_k)|| - (p-n_1-2)\geq {{p^2}\over 4},$
or equivalently, $$ { {4 -2(p-n_1)+\sum _{i=1}^{k}  n_i(p-n_i)}}\geq {{p^2}\over 2}.$$ 

By hypothesis, for every integer $i=1,2,\ldots, k$,  $n_1\geq n_i$, and as a consequence $p-n_i\geq p-n_1$. Then, $\sum _{i=1}^{k}  n_i(p-n_i)\geq (p-n_1) \sum _{i=1}^{k}n_i=p(p-n_1)$ and it is sufficient to show %{\color{blue} $$4 -2(p-n_1)+\sum _{i=1}^{k}  n_i(p-n_i)\geq 4 -2(p-n_1)+p(p-n_1)=4 + (p-2)(p-n_1)$$ } 
$$ { 4 + (p-2)(p-n_1)}\geq {{p^2}\over 2}.$$

 \noindent {\bf Case 1)} If $n_1\leq {{p-2}\over 2}$, then $p-n_1\geq {{p+2}\over 2}$ and we conclude
 $$ { 4 + (p-2)(p-n_1)} \mbox{    } \geq \mbox{    } { 4 + (p-2)\Big({{p+2}\over 2}\Big)} \mbox{    }=\mbox{    } {2+{{p^2}\over 2}}.$$
 
 \noindent {\bf Case 2)} If ${{p-2}\over 2} < n_1$, imposing the necessary condition $n_1\leq {p\over2}$ leaves only two possibilities: $p$ is even and $n_1={p\over 2}$, or $p$ is odd and $n_1={{p-1}\over 2}$. In these cases we write  $${4 -2(p-n_1)+\sum _{i=1}^{k}  n_i(p-n_i)}\mbox{    } = \mbox{    } {4 -2(p-n_1)+n_1(p-n_1)+\sum _{i=2}^{k}  n_i(p-n_i)}.$$ 
 By hypothesis, again, for every integer $i$, $i=2,\ldots, k$, $n_2\geq n_i$, and this implies $p-n_i\leq p-n_2$. Also, note that $\sum _{i=2}^{k} n_i=p-n_1$, and as a consequence, $\sum _{i=2}^{k}  n_i(p-n_i)\geq (p-n_2)\sum _{i=2}^{k}  n_i=(p-n_2)(p-n_1)$. Then, it is sufficient to prove
 $$ { 4 +(n_1-2)(p-n_1)+(p-n_2)(p-n_1)} \mbox{    } \geq \mbox{    } {{p^2}\over 2}.$$
 
Suppose $n_2\leq n_1-2$, so that $p-n_2\geq p-(n_1-2)$. Then,  $$ { 4 +(n_1-2)(p-n_1)+(p-n_2)(p-n_1)} \mbox{    } \geq  \mbox{    }  4+p(p-n_1)$$ and it is sufficient to check that when $n_1={p\over 2}$ ($p$ even) or $n_1={{p-1}\over 2}$ ($p$ odd), then $4+p(p-n_1)\geq {{p^2}\over 2}$, a straightforward verification.

%{\color{blue}    If $n_1={p\over 2}$, then $4+p(p-n_1)=4+p\big({p\over 2}\big)=4+{{p^2}\over 2}$\\ \indent If  $n_1={{p-1}\over 2}$, then $4+p(p-n_1)=4+p\big({{p+1}\over 2}\big)=4+{{p^2+p}\over 2}$}
\smallskip
 If $n_2> n_1-2$, since $n_2\leq n_1$, either $n_1=n_2$ or $n_2=n_1-1$ and the only possibilities are: 
\begin{itemize}
\item $p$ even, $k=3$, $n_1={p\over 2}$, $n_2= {p\over 2}-1$ and $n_3=1$ %{\color{blue}  $$ 4-2p+2n_1+n_1(p-n_1)+n_2(p-n_2)+n_3(p-n_3)=$$ $$4-2p+p+\Big({p\over 2}\Big)^2+\Big({p\over 2}-1\Big)\Big({p\over 2}+1\Big)+p-1 = 4+\Big({p\over 2}\Big)^2+\Big({p\over 2}\Big)^2-1-1= 2+{{p^2}\over 2}$$}
\item $p$ odd, $k=3$, $n_1={{p-1}\over 2}$,  $n_2={{p-1}\over 2}$ and $n_3=1$ %{\color{blue}  $$ 4-2p+2n_1+n_1(p-n_1)+n_2(p-n_2)+n_3(p-n_3)=$$ $$4-2p+p-1+2\Big({{p-1}\over 2}\Big)\Big({{p+1}\over 2}\Big)+p-1 = 2+{{(p-1)(p+1)}\over 2}= 2+{{p^2-1}\over 2} = {3\over 2}+{{p^2}\over 2}$$}
\item $p$ odd, $k=3$, $n_1={{p-1}\over 2}$,  $n_2={{p-3}\over 2}$ and $n_3=2$% {\color{blue} $$4 -2p+2n_1+n_1(p-n_1)+n_2(p-n_2)+n_3(p-n_3)=$$ $$ 4-2p+p-1+\Big({{p-1}\over 2} \Big)\Big({{p+1}\over 2}\Big)+\Big({{p-3}\over 2} \Big)\Big({{p+3}\over 2}\Big)+ 2(p-2)= {p-1+{{p^2-5}\over 2}}$$ $ > {{p^2}\over 2}+p-1- 3$ and for this case to happen, since $n_1\geq n_2\geq n_3=2$ and $p$ is odd, $p\geq 7$ and $p-4>0$ }
\item $p$ odd, $k=4$, $n_1={{p-1}\over 2}$,  $n_2={{p-3}\over 2}$, $n_3=1$ and $n_4=1$ %{\color{blue} $$ 4-2p+2n_1+n_1(p-n_1)+n_2(p-n_2)+n_3(p-n_3)=$$  $$4 - 2p +p-1+\Big({{p-1}\over 2} \Big)\Big({{p+1}\over 2}\Big)+\Big({{p-3}\over 2} \Big)\Big({{p+3}\over 2}\Big)+2(p-1)= {p+1+ {{p^2-1}\over 4} +{{p^2-9}\over 4} }=$$
%$ > {{p^2}\over 2}+p+1- 3> {{p^2}\over 2}$ }
\end{itemize}
In these cases, it is straightforward to verify ${ {4 -2(p-n_1)+\sum _{i=1}^{k}  n_i(p-n_i)}}\geq {{p^2}\over 2}.$
\epf

\bigskip

We close by noting that in \cite{DKPE17}, DeBiasio et al. considered minimum degree conditions corresponding to Corollary 12 for $k$-partite graphs that are {\it fair}.  They used a parameter $\lambda$ and gave asymptotic results on $\lambda$-fair graphs with $\lambda\geq 2$. Our necessary condition, $n_1\leq {p\over 2}$ is equivalent to $\lambda \geq 2$ but $G$ does not need to be $\lambda$-fair for our results to hold. 

\section*{Acknowledgments}
We thank IPFW for its hospitality and Jay Bagga, Lowell Beineke and Marc Lipman for their comments at the beginning of this project. We thank the MAA - Tensor Program for Women and Mathematics for facilitating a research visit during this project. We also thank the anonymous referee who provided insightful comments, brought to our attention the work in \cite{C18} and suggested the study of chorded pancyclicity.


\smallskip

\begin{thebibliography}{20}

\bibitem{A09} J. Adamus. Edge condition for hamiltonicity in balanced tripartite graphs. {\it Opuscula Math.} {\bf 29}, 337--343 (2009) 

\bibitem{B71} J.A. Bondy. Pancyclic graphs I. {\it J. Combin.Theory Ser. B} {\bf 11}, 80--84 (1971) 

\bibitem{CFGJL95} G. Chen, R. Faudree, R. Gould, M.S. Jacobson, L. Lesniak. Hamiltonicity of balanced $k$-partite graphs. {\it Graphs Combin.} {\bf 11}, 221--231 (1995) 

\bibitem{C18} G. Chen, R. J. Gould, X. Gu, A. Saito. Cycles with a chord in dense graphs. {\it Discrete Math.} {\bf 341}, 2131--2141 (2018)  %%% added 

\bibitem{CJ97}  G. Chen, M.S. Jacobson. Degree sum conditions for hamiltonicity on $k$-partite graphs. {\it Graphs Combin.} {\bf 13}, 325--343 (1997) 
 
\bibitem{C17} M. Cream, R. J. Gould, K. Hirohata. A note on extending Bondy's meta-conjecture. {\it Australas. J. Combin.} {\bf 67(3)}, 463--469 (2017)    %%% added 

\bibitem{DKPE17} L. DeBiasio, R. A. Krueger, D. Pritikin, E. Thompson. Hamiltonian cycles in fair $k$-partite graphs. arXiv:1707.07633v2 [math.CO] (2017)

\bibitem{D52} G.A. Dirac. Some theorems on abstract graphs. {\it Proc. Lond. Math. Soc. } {\bf (3)2}, 69--81 (1952) 

\bibitem{ES88} R.C. Entringer,  E. Schmeichel. Edge conditions and cycle structure in bipartite graphs. {\it Ars Combin.} {\bf 26}, 229--232 (1988) 

\bibitem{GJ79} M. R. Garey, D. S. Johnson. {\it Computers and Intractability: A Guide to the Theory of NP-Completeness}, W.H. Freeman \& Co. New York, NY (1979)

\bibitem{MM63} J. Moon, L. Moser.  On hamiltonian bipartite graphs. {\it Israel J. Math.} {\bf 1}, 163--165 (1963)

\bibitem{O60} O. Ore. Note on hamilton circuits. {\it Amer. Math. Monthly} {\bf 67(1)}, 55-55 (1960)

\bibitem{P62} L. P{\'o}sa. A theorem concerning hamiltonian lines. {\it Alkalmaz. Mat. Lapok.} {\bf 7},  225--226 (1962) 

\end{thebibliography}
\end{document}